\documentclass[a4paper,1àpt]{article}

\usepackage{euscript}
\usepackage{amsfonts}
\usepackage{amsmath,amscd, amsthm, amssymb}
\usepackage[english]{babel}
\usepackage[latin1]{inputenc}
\usepackage{graphicx}
\usepackage[arrow, matrix, curve]{xy}
\usepackage{dsfont}\let\mathbb\mathds
\usepackage{rotating}

\newtheorem{thm}{Theoreme}
\newtheorem{prop}[thm]{Proposition}

\newtheorem{lem}[thm]{Lemma}

\newtheorem{Def}[thm]{Definition}

\def\cartesien{%
    \ar@{-}[]+R+<6pt,-1pt>;[]+RD+<6pt,-6pt>%
    \ar@{-}[]+D+<1pt,-6pt>;[]+RD+<6pt,-6pt>%
  }
\newcommand{\immouv}[1][r]
   {\ar@{}[#1] |*[o][F]{\hbox{%
         \vrule width 1.5mm height 0pt depth 0pt%
         \vrule width 0pt height .75mm depth .75mm%
         }}
     \ar@{^{(}->}[#1]}
\newcommand{\demde}[1]{\begin{proof} de #1}
\newcommand{\dem}{\begin{proof}}
\newcommand{\cqfd}{\end{proof}}
\let\cat\mathfrak 
  \newcommand{\UN}[4][r]{%
    \ar@/^1pc/[#1]^{#2}_*=<0.3pt>{}="HAUT"
    \ar@/_1pc/[#1]_{#3}^*=<0.3pt>{}="BAS"
    \ar @{=>} "HAUT";"BAS" ^{#4}
  }
\newcommand{\DEUX}[6][r]{
    \ar@/^2pc/[#1]^{#2}_*=<0.3pt>{}="HAUT"
    \ar@{}    [#1]     ^*=<0.3pt>{}="MILIEUHAUT"
                       _*=<0.3pt>{}="MILIEUBAS"
    \ar[#1]_(0.3){#3}                  
    \ar@/_2pc/[#1]_{#4}^*=<0.3pt>{}="BAS"
    \ar @{=>} "HAUT";"MILIEUHAUT" ^{#5}
    \ar @{=>} "MILIEUBAS";"BAS" ^{#6}
  }   
 \newcommand{\eq}[1][r]
   {\ar@<-3pt>@{-}[#1]
    \ar@<-1pt>@{}[#1]|<{}="gauche"
    \ar@<+0pt>@{}[#1]|-{}="milieu"
    \ar@<+1pt>@{}[#1]|>{}="droite"
    \ar@/^2pt/@{-}"gauche";"milieu"
    \ar@/_2pt/@{-}"milieu";"droite"}
 \newcommand{\incl}[1][r]
  {\ar@<-0.2pc>@{^(-}[#1] \ar@<+0.2pc>@{-}[#1]}

\newcommand{\Ac}{\mathcal{A}}
\newcommand{\Bc}{\mathcal{B}}
\newcommand{\Cc}{\mathcal{C}}

\newcommand{\Ic}{\mathcal{I}}
\newcommand{\Jc}{\mathcal{J}}

\newcommand{\Tc}{\mathcal{T}}

\newcommand{\aaa}{\mathfrak{a}}

\newcommand{\bbb}{\mathfrak{b}}
\newcommand{\ccc}{\mathfrak{c}}

\def\dar[#1]{\ar@<2pt>[#1]\ar@<-2pt>[#1]}
  \entrymodifiers={!!<0pt,0.7ex>+}

\title{Interchange of filtered $2$-colimits and finite $2$-limits}

\author{Delphine Dupont}

\begin{document}
\date{}
\maketitle
\begin{abstract}
In this paper we go into the study of $2$-limit and $2$-colimit in the $2$-category $\Cc\Ac\Tc$ the category of small categories. In particular we show the commutation of filtered $2$-colimits and finite $2$-limits. It is a generalization of a classical result in category theory (see for example \cite{Macl}).
\end{abstract}

\section{Introduction}
Recently the 2-category theory has developed rapidly. It gives a very useful language in many field of mathematics. One of the main notions is the $2$-limits and $2$-colimits. For example the stalk of a stack is a $2$-colimit. In this paper we genralize a classical result in category theory.   

In a first part we recall the definition of  $2$-limit and $2$-colimit in $\Cc\Ac\Tc$, the $2$-category of small categories. We give an explicit description of a filtered $2$-colimit and of a $2$-limit. But we do not express the morphisms of these categories in a classical way. We give them as elements of a colimit or a limit. 

In a second part we prove the interchange of filtered $2$-colimits and finite $2$-limits. To prove this theorem we use the explicit expression of  the categories  $\begin{displaystyle}
2 \varinjlim_{i \in \Ic}2\varprojlim_{j \in \Jc}\aaa(i,j) \end{displaystyle}$ and  $\begin{displaystyle}  2\varprojlim_{j \in \Jc}2\varinjlim_{i \in \Ic}\aaa(i,j)
 \end{displaystyle}$.

We do not recall the definitions of a $2$-category, a $2$-functor (some times called pseudo functor in the literature), a $2$-natural transformation of functors and a $2$-modification. The reader can find them in the chapter 7 of  \cite{Macl}, the appendix of the paper \cite{Ingo1} and the paper \cite{St}. \\

Let us recall the definition of a filtered category : 
\begin{Def}
A category $\Ic$ is filtered if it satisfies the conditions (i)-(iii) below.
\begin{itemize}
\item[(i)] $\Ic$ is non empty,
\item[(ii)] for any $i$ and $j$ in $\Ic$, there exist $k \in \Ic$ and morphisms $i\rightarrow k$, $j\rightarrow k$,
\item[(iii)] for any parallel morphisms $s,s'  :\raisebox{.7ex}{\xymatrix{ i \dar[r] & j }}$ there exists a morphism $h : j\rightarrow k$ such that $h\circ s= h\circ s'$.
\end{itemize}
\end{Def}

\section{$2$-colimits and $2$-limits}
Let us first recall the definitions of a $2$-limit and a $2$-colimit. This part is inspired by the appendix of  \cite{Ingo1}. We cite \cite{St} as a classical reference.

\begin{Def}
Let $\Ic$ be a small category, and $\cat b: \Ic\rightarrow \Cc \Ac \Tc $ a $2$-functor
 The system $\cat b$ admits a $2$-colimit if and only if there exist :
\begin{itemize}
\item a category $2\varinjlim_{i\in \Ic}\cat b(i)$ and
\item a $2$-natural transformation $\sigma :\cat b \rightarrow 2\varinjlim_{i\in \Ic}\cat b(i )$ 
\end{itemize}
such that for all category $\Cc$ the functor :
$$( \sigma \circ) : \text{\underline{Hom}}_{\Cc}(2\varinjlim_{i\in \Ic}\cat b(i), \Cc) \rightarrow \text{\underline{Hom}}(\cat b, \Cc)$$
is an equivalence of categories. \\
We say that a $2$-colimit has the strong factorisation property if $( \sigma\circ)$ is an isomorphism of categories.
\end{Def}
More concretely, a $2$-functor $\cat b : \Ic\rightarrow \Cc\Ac\Tc$ admits a $2$-colimit if and only if there exist :
\begin{itemize}
\item a category $2\varinjlim_{i \in \Ic}\cat b(i)$,
\item functors $\sigma_{i} : \cat b(i)\rightarrow 2\varinjlim_{i \in \Ic}\cat b(i)$ for any $i \in \Ic$,
\item and a natural equivalence $\Theta_{s}^\sigma : \sigma_{i}\buildrel\sim\over\rightarrow \sigma_{j}\cat b(s)$ for any morphism $s: i\rightarrow j$ of $\Ic$ visualized by :
$$
\xymatrix @!0 @C=0.7cm @R=0.7cm {
&&&2\varinjlim_{i\in \Ic}\cat b(i)  \ar[ddddlll]_{\sigma_{i}}  \ar[ddddrrr]^{\sigma_{j}}\\
~\\
&&&&~\\
&&\ar@{=>}[urr]_{\theta_{s}}^\sim&&\\
\cat b(i) \ar[rrrrrr]_{\cat b(s)}&&&&&&\cat b(j)\\
}$$
such that $\sigma :\cat b \rightarrow \varinjlim_{i\in \Ic}\cat b(i)$ is a $2$-natural transformation.
\end{itemize}
Moreover these data satisfy condition assuring that $(\circ \sigma)$ is an equivalence. The first one translate the fact that $(\circ \sigma)$ is essentially surjective and the second one that $(\circ \sigma)$ is fully faithfull.
\begin{itemize}
\item For any category $\Cc$, any morphism of functors $\rho : \cat b\rightarrow \Cc$, there exist a functor :
$$F: 2\varinjlim_{i\in \Ic}\cat b(i) \longrightarrow \Cc$$
and an isomorphism $\varphi^F:\rho \rightarrow F\sigma$, which is a modification given, for all $i \in \Ic$, by a natural equivalence   $\varphi_{i}^F :\rho_{i} \rightarrow F\circ \sigma_{i}$. This may be visualized by :

$$
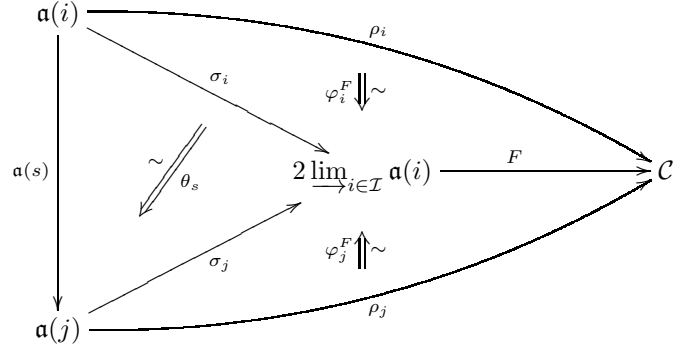
\xymatrix @!0 @C=1cm @R=0.7cm {
\cat a(i)\ar[dddddd]_{\cat a(s)} \ar[dddrrrr]^{\sigma_{i}} \ar@/^1.5pc/[rrrrrrrrddd]^{\rho_{i}}\\
&&&& \ar@{=>}[d]^\sim_{\varphi_{i}^F}\\
&& \ar@{=>}[ddl]^{\theta_{s}}_\sim&&\\
&&&&2\varinjlim_{i\in \Ic}\cat a(i) \ar[rrrr]^F&&&&\Cc\\
&&&&&&\\
&&&&\ar@{=>}[u]^{\varphi_{j}^F}_{\sim}\\
\cat a(j)\ar[uuurrrr]_{\sigma_{j}} \ar@/_1.5pc/[uuurrrrrrrr]_{\rho_{j}}
}$$
The compatibility condition is given by :

$$\big(F\bullet\theta_{s}^\sigma\big)\circ \varphi_{i}^F=\big(\varphi_{j}^F\bullet\cat a(s)\big)\circ \theta_{s}^\rho$$
The pair $(F, \varphi^F)$ is called a lax factorization of the system $\rho$.
\item[$\bullet$]  Let $\rho$ and $\rho'$ be two $2$-natural transformations and $\lambda : \rho\rightarrow \rho'$ a modification. Then for any lax factorization $F: 2\varinjlim_{i\in \Ic}\cat b(i) \longrightarrow \Cc$ of $\rho$ and $G : 2\varinjlim_{i\in \Ic}\cat b(i) \longrightarrow \Cc$ of $\rho'$ there exists a unique natural transformation $\Lambda : F\rightarrow G$ such that :
$$\varphi_{i}^G\circ \lambda_{i}=\big(\Lambda_{i}\bullet\sigma_{i}\big)\circ \varphi_{i}^{F}$$
\end{itemize}
If the $2$-colimit has the strong factorization property then there exists a
unique factorization such that $\varphi^F$ is the identity.\\\\
A contravariant $2$-functor $\cat c :\Ic \rightarrow \Cc \Ac\Tc$ is a $2$-functor $\cat c :\Ic^{op}\rightarrow \Cc\Ac\Tc$. A $2$-limit is construction dual of the $2$-colimit. Hence, we have the following definition :
\begin{Def}
Let $\Ic$ be a small category and $\cat c: \Ic\rightarrow \Cc \Ac \Tc $ a contravariant $2$-functor. The system $\cat c$ admits a $2$-limit if and only if there exist :
\begin{itemize}
\item a category $2\varprojlim_{i\in \Ic}\cat c(i)$ and
\item a $2$-natural transformation $\sigma : 2\varprojlim_{i\in \Ic}\cat b(i )\rightarrow \cat c$ 
\end{itemize}
such that for all category $\Cc$ the functor :
$$( \circ\sigma ) : \text{\underline{Hom}}(\cat c, \Cc)\rightarrow   \text{\underline{Hom}}_{\Cc}(2\varprojlim_{i\in \Ic}\cat b(i), \Cc)$$
is an equivalence of categories. \\
We say that a $2$-limit has the strong factorisation property if $(\circ \sigma)$ is an isomorphism of categories.
\end{Def}

It is well known that the $2$-category $\Cc \Ac \Tc$ is complete and co-complete. The proof of this result consist an explicit definition of the $2$-limit and $2$-colimit. We are going to recall these definitions, but morphisms of these categories won't be given in a classical way. Usually they are given as classes of morphisms between objects of $\bbb(i)$ which satisfy some conditions.  Here we are going to express them in term of limit and colimit.

Let $\bbb: \Ic \rightarrow \Cc\Ac \Tc$ be a $2$-functor.
Let us give some useful notations.\\
Let $i, i' \in \Ic$ and $\Ic_{ii'}$ be the category defined by:
\begin{itemize}
\item the objects of $\Ic_{i,i'}$ are :
$$\Big\{ (i'', s,s')\mid i''\in Ob\Ic, s \in Hom_{\Ic}(i,i''), s'\in Hom_{\Ic}(i',i'')\Big\}$$
visualized by :
$$
\xymatrix @!0 @C=1.5cm @R=0.4cm {
i\ar[dr]^s\\
& i''\\
i'\ar[ur]_{s'}
}$$
\item the morphisms of $\Ic_{i,i'}$ from $(i_{1}'',s_{1},s_{1}')$ to $(i_{2}'',s_{2}, s_{2}' )$ are : 
$$Hom_{\Ic_{ii'}}\big((i_{1}'',s_{1},s_{1}'), (i_{2}', s_{2}, s_{2}')\big)=\{t\in Hom_{\Ic}(i''_{1},i''_{2})\mid t\circ s_{1}=s_{2}, t\circ s_{1}'=s_{2}'\}  $$
visualized by :
$$
\xymatrix @!0 @C=1.5cm @R=0.9cm {
i\ar[dr]^{s_{1}} \ar@/^/[rrd]^{s_{2}}\\
& i''_{1}\ar[r]^t& i''_{2} \\
i'\ar[ur]_{s_{1}'} \ar@/_/[rru]_{s_{2}'}
}$$
\end{itemize}
\begin{lem}
If the category $\Ic$ is filtered then the category $\Ic_{ii'}$ is filtered.
\end{lem}
\dem
The proof is straight-foward.
\cqfd

\begin{prop}\label{colimit}
The category $\Bc$ defined below is a $2$-colimit of $\bbb$ satisfying the strong factorization property.
\begin{itemize}
\item Objects of $\Bc$ are pairs $(i,X)$ where $i \in \Ic$ and $X \in \cat b(i)$,
\item let $(i,X)$ and $(i',Y)$ be two objects of $\Bc$, the morphisms from $(i,X)$ to $(i', Y)$ are the elements of the colimit :
$$\begin{displaystyle}Hom_{\Bc}\big(X,Y\big)=\varinjlim_{(i'',s,s')\in \Ic_{ii'}}Hom_{\bbb(i'')}\big(\bbb(s)X,\bbb(s')Y\big)\end{displaystyle}$$
where if $t \in Hom_{\Ic_{ii'}}\big((i''_{1},s_{1}, s'_{1}), (i_{2}'', s_{2}, s_{2}')\big)$, its image in the inductive system is the following morphism :
$$\begin{array}{cccc}
 Hom_{\cat b(i_{1}'')}\big(\cat b(s_{1})(X), \cat b(s_{1}')(Y)\big)& \longrightarrow & Hom_{\cat b(i_{2}'')}\big(\cat b(s_{2})(X), \cat b(s_{2}')(Y)\big)\\
 h& \longmapsto & b_{t,s_{1}'}^{-1}\circ \cat b(t)h\circ b_{t, s_{1}}
\end{array}$$
where $b_{t,s_{1}}$ is the isomorphism given by the $2$-functor $\cat b$ :
$$b_{t,s_{1}} : \cat b (t \circ s_{1}) \buildrel\sim\over\longrightarrow \cat b(t)\circ \cat b(s_{1})$$

\end{itemize}
\end{prop}
\dem
See for example \cite{Ingo1}.
\cqfd
If $\Ic_{ii'}$ is filtered, the morphisms from $(i, X)$ to $(i', Y)$ are the elements of the quotient : 
$$\coprod_{i'' \in \Ic_{ii'}}Hom_{\cat b(i'')}\big(\bbb(s)X, \bbb(s')Y\big)\Big/ \sim$$
Where, given $h_{1}\in  Hom_{\cat b(i_{1}'')}\big(\cat b(s_{1})X, \cat b(s_{1}')Y\big)$ and $h_{2}\in  Hom_{\cat b(i_{2}'')}\big(\cat b(s_{2})X, \cat b(s_{2}')Y\big)$, we write $h_{1}\sim h_{2}$ if and only if there exist $(i_{3}, s_{3}, s_{3}')\in \Ic_{ii'}$, $t_{13} \in Hom_{}\big((i_{1},s_{1},s_{1}'),(i_{3}, s_{3}, s_{3}')\big)$ and $t_{23} \in Hom_{}\big((i_{2},s_{2},s_{2}'),(i_{3}, s_{3}, s_{3}')\big)$ such that :
$$b_{t_{13},s_{1}'}^{-1}\circ\bbb(t_{13})h_{1}\circ b_{t_{13}, s_{1}}=b_{t_{23},s_{2}'}^{-1}\circ\bbb(t_{23})h_{2}\circ b_{t_{23}, s_{2}}$$
Similarly we are going to give an explicit construction of a $2$-limit. As before morphisms of this category will be given by a limit. Let $\ccc$ be a contravariant $2$-functor :
$$\ccc : \Jc^{op} \longrightarrow \Cc \Ac \Tc$$
\begin{prop}\label{limit}
The  category $\Cc$ define below, is a $2$-limit of $\ccc$ satisfying the strong factorization property.
\begin{itemize}
\item The objects of $\Cc$ are pairs $(X, \vartheta^X)$ where :
\begin{itemize}
\item $X=\{X_{j}\}_{j \in \Jc}$, for $X_{j} \in Ob(\ccc(j))$,
\item $\vartheta^X=\{ \vartheta^X_{t}\}_{t \in Hom_{\Jc}}$ where, for  $t\in Hom_{\Jc}(j,j')$, $\vartheta_{t}^X$ is an isomorphism :
$$\vartheta_{t}^X : X_{j}  \buildrel\sim\over\longrightarrow \ccc(t)X_{j'}$$
\end{itemize}
satisfying the following conditions :
\begin{itemize}
\item[A)] for any $j \in \Jc$ we have $Id_{j}= c_{j}(X_{j})\circ\vartheta_{Id_{j}}^X$, where $c_{j}$ is the morphism given by the $2$-functor $\cat c$ :
$$c_{j}:\ccc(Id_{j})\buildrel\sim\over\longrightarrow Id_{\ccc(j)}$$
\item[B)] for any two composable morphisms $t : j\rightarrow j'$ and $t': j'\rightarrow j''$ the following equation holds :
$$\ccc(t)( \vartheta_{t'}^X)\circ \vartheta_{t}^X=c_{t,t'}(X_{j''})\circ \vartheta_{t\circ t'}^X$$
where $c_{t,t'}$ is the isomorphism given by the $2$-functor $\cat c$ :
$$c_{t,t'}:\cat c(t\circ t') \buildrel\sim\over\longrightarrow \cat  c(t) \circ \cat c(t')$$
\end{itemize}

\item Let $(X, \vartheta^{X})$ and $(Y, \vartheta^{Y})$ be two objects of $\Cc$, morphisms from $(X, \vartheta^X)$ to $(Y, \vartheta^Y)$ are elements of the limit : 
$$\begin{displaystyle}Hom_{\Cc}(X,Y)=\varprojlim_{j \in \Jc}Hom_{\ccc(j)}\big(X_{j}, Y_{j}\big)\end{displaystyle}$$
where if $t \in Hom_{\Jc}(j,j')$, its image in the projective system is the following morphism :
$$\begin{array}{cccc}
 Hom_{\cat c(j)}(X_{j}, Y_{j}) &\longrightarrow & Hom_{\cat c(j')}(X_{j'}, Y_{j'})\\
  h & \longmapsto & (\vartheta_{t}^Y)^{-1}\circ\ccc(t)(h)\circ \vartheta_{t}^X
\end{array}$$
\end{itemize}
\end{prop}
This means that a morphism between two objects is the datum of $\{h_{j}\}_{j \in \Jc}$, where $h_{j} \in Hom_{\cat c(j)}(X_{j}, Y_{j})$ satisfies the equality :
$$ h_{j}=(\vartheta_{t}^Y)^{-1}\circ\ccc(t)(h_{j'})\circ \vartheta_{t}^X$$
for all $t :j'\rightarrow j$ morphisms of $\Jc$.

\section{Interchange of filtered $2$-colimits and finite $2$-limits}
We are going to show that filtered $2$-colimits commute with finite $2$-limits. Let us give first a precise meaning to this sentence.\\
Let $\Ic$ be a filtered category, $\Jc$ a finite category and $\aaa$ a $2$-functor : 
$$\aaa : \Ic\times \Jc^{op} \longrightarrow \Cc \Ac \Tc$$
Let us denote $2\bold{\mathfrak{F}}(\Cc)$ the $2$-category of $2$-functors going from the category $\Cc$ to the $2$-category $\Cc \Ac \Tc$. We have the following proposition, for a proof see for example \cite{Ingo1} :
\begin{prop}
The correspondence :
$$
\begin{array}{cccc}
2\mathfrak{F}(\Cc) & \longrightarrow & \Cc \Ac \Tc\\
\bbb & \longmapsto  & \begin{displaystyle}2\lim_{\substack{ \longrightarrow \\ i \in \Ic }} \bbb(i)\end{displaystyle}
\end{array}
$$
can be extended to a $2$-functor between $2$-categories. A similar statement holds for $2$-limits.
\end{prop}
\noindent
Now let us consider, the natural $2$-functor :
$$\begin{array}{ccc}
\Jc & \longrightarrow & 2\mathfrak{F}(\Ic)\\
j & \longmapsto & \aaa(\cdot, j).

\end{array}$$
 The composition of this $2$-functor and the one defined in the proposition gives a functor from $\Jc$ to $\Cc\Ac\Tc$. As $\Cc\Ac\Tc$ is complete we can consider its limit. Let us denote 
 $$\begin{displaystyle} 2\varprojlim_{j \in \Jc}2\varinjlim_{i \in \Ic}\aaa(i,j)\end{displaystyle}$$ this limit. We define in the same way the $2$-colimit $\begin{displaystyle}2\varinjlim_{i \in \Ic}2\varprojlim_{j \in \Jc}\aaa(i,j)\end{displaystyle}$. \\
Remark\\
The composition, for all $i\in\Ic$ and $j \in \Jc$, of the functor define by the $2$-colimit and the $2$-limit : 
$$ i \begin{displaystyle}2\varprojlim_{j' \in \Jc}\aaa(i,j') \longrightarrow \cat a(i,j)\longrightarrow 2\varinjlim_{i' \in \Ic}\aaa(i',j)i \end{displaystyle}$$
defines a functor :
$$\Psi :\begin{displaystyle}
2 \varinjlim_{i \in \Ic}2\varprojlim_{j \in \Jc}\aaa(i,j) \longrightarrow  2\varprojlim_{j \in \Jc}2\varinjlim_{i \in \Ic}\aaa(i,j)
 \end{displaystyle}$$
 \begin{thm}
 The natural functor
  $$\Psi : \begin{displaystyle}
2 \varinjlim_{i \in \Ic}2\varprojlim_{j \in \Jc}\aaa(i,j) \longrightarrow  2\varprojlim_{j \in \Jc}2\varinjlim_{i \in \Ic}\aaa(i,j)
 \end{displaystyle}$$
 is an equivalence of categories. 
 \end{thm}
\dem
The fact that $\Psi$ is fully faithful comes directly from the expression of the morphisms of these two categories in terms of limits and colimits. \\
In detail, using the two propositions \ref{colimit} and \ref{limit}, we can give an explicit construction of the categories $\begin{displaystyle} 2\varinjlim_{i \in \Ic}2\varprojlim_{j \in \Jc}\aaa(i,j) \end{displaystyle}$ and $\begin{displaystyle}2\varprojlim_{j \in \Jc}2\varinjlim_{i \in \Ic}\aaa(i,j) \end{displaystyle}$.

The category
 $\begin{displaystyle} \varinjlim_{i \in \Ic}\varprojlim_{j \in \Jc}\aaa(i,j) \end{displaystyle}$ is defined as follows :
  \begin{itemize}
 \item[$\bullet$] its objects  are the triples 
 $\big(i,X, \vartheta^X \big)$ where 
 \begin{itemize}
 \item $i $ is an object of  $\Ic$, 
 \item $X=\{X_{j}\}_{j \in \Jc}$ for $X_{j}$ is an object of $ \aaa(i,j)$ 
 \item $\vartheta^X=\{\vartheta_{t}^X\}_{t \in Hom\Jc}$, where, for  $t: j' \rightarrow j$ an isomorphism of $\Jc$, $\vartheta_{t}^X$ is a morphism :
 $$\vartheta_{t}^X : X_{j} \longrightarrow \aaa(Id_{i}, t)X_{j'}$$
 \end{itemize}
 verifying the two following conditions : 
\begin{itemize}
\item for any $j \in Obj \Jc$ we have $a_{ij}(X_{ij})\circ \vartheta_{Id_{i}}^X= Id_{X_{i}}$, where $a_{ij}$ is the isomorphism $a_{ij} : \aaa(Id_{i}, Id_{j}) \buildrel\sim\over\rightarrow Id_{\aaa(i,j)}$ given by the $2$-functor $\aaa$,
\item for any two composable morphisms $t : j \rightarrow j'$ and $t' : j' \rightarrow j''$ the following equation holds :
$$\aaa(Id_{i}, t)(\vartheta_{t'}^X) \circ \vartheta_{t}^X= a_{(Id_{i}, t), (Id_{i}, t')}(X_{ij''})\circ \vartheta_{t'\circ t}$$
where $a_{(Id_{i}, t), (Id_{i}, t')} : \aaa(Id_{i},t' \circ t)\buildrel\sim\over\rightarrow \aaa(Id_{i}, t')\aaa(Id_{i}, t)$ is the morphism given by the $2$-functor $\aaa$.
\end{itemize} 
\item[$\bullet$] the set of  morphisms from $\big(i, X,\vartheta^X\big)$ to $\big(i',Y, \vartheta^Y\big)$ is given by the limit :
$$Hom_{}(X,Y):=\begin{displaystyle}\varinjlim_{i''\in \Ic_{ii'}}\varprojlim_{j \in \Jc}Hom_{\cat a(i'',j)}\big(\cat a(s, Id_{j})X_{j}, \cat a(s', Id_{j})Y_{j}\big)\end{displaystyle}.$$
\end{itemize}
The category
 $\begin{displaystyle}2\varprojlim_{j \in \Jc}2 \varinjlim_{i \in \Ic}\aaa(i,j) \end{displaystyle}$ is defined as follows :
\begin{itemize}
\item[$\bullet$] the objects are pairs $(\mathcal{X}, \theta^\mathcal{X})$, where 
\begin{itemize}
\item $\mathcal{X}=\{(i_{j}, X_{i_{j}})\}_{j \in \Jc}$ with $i_{j} \in \Ic$ and $X_{i_{j}} \in \cat a(i_{j},j)$,
\item $\theta^\mathcal{X}=\{[\vartheta_{t}]\}_{t \in Hom\Jc}$ where, for $t \in Hom_{\Jc}(j,j')$, $[\vartheta_{t}^\mathcal{X}]$ belongs to the quotient :
$$\coprod_{i'' \in \Ic_{i_{j}i_{j'}}}Hom_{\cat a(i'', j)}\big(\aaa(s, Id_{j})X_{i_{j}}, \aaa(s',t)X_{i_{j}}\big)\Big/ \sim$$
and where $[\vartheta_{t}^\mathcal{X}]$ 
satisfies the following equalities :
\begin{equation}[a_{ij}(X_{ij})\circ \vartheta_{Id_{i}}^X]= [Id_{X_{i}}]\end{equation}
\begin{equation}[\aaa(Id_{i}, t)(\vartheta_{t'}^X) \circ \vartheta_{t}^X]=[ a_{(Id_{i}, t), (Id_{i}, t')}(X_{ij''})\circ \vartheta_{t'\circ t}]\end{equation}
and for any two composable morphisms $t : j \rightarrow j'$ and $t' : j' \rightarrow j''$.
\end{itemize}
\item[$\bullet$] the set of morphisms from $(\mathcal{X}, \theta^\mathcal{X})$ to $(\mathcal{Y}, \theta^\mathcal{Y})$  is given by the limit :
$$\begin{displaystyle}
\varprojlim_{j \in \Jc}\varinjlim_{i'' \in \Ic_{i_{j}i_{j}'}}  Hom_{\cat a(i'',j)}\big(\cat a(s, Id_{j})X_{j
}, \cat a(s', Id_{j})Y_{j}\big)
\end{displaystyle}$$
\end{itemize} 
The natural functor between the $2$-limits is  :
$$\begin{array}{ccccc}
\varepsilon :&\begin{displaystyle} 2\varinjlim_{i \in \Ic}2\varprojlim_{j \in \Jc}\aaa(i,j)\end{displaystyle}& \longrightarrow &\begin{displaystyle}2\varprojlim_{j \in \Jc}2 \varinjlim_{i \in \Ic}\aaa(i,j)\end{displaystyle}\\
&\big(i, \{X_{j}\}, \{\vartheta^X_{t}\}\big) & \longmapsto & \big(\{i, X_{j} \}, \{[\vartheta_{t}^X]\}\big)\\
\end{array}$$
Moreover, if $\mathcal{X}=(i,X,\vartheta^X)$ and $\mathcal{Y}=(i',Y, \vartheta^Y)$ are two objects of $ \begin{displaystyle} 2\varinjlim_{i \in \Ic}2\varprojlim_{j \in \Jc}\aaa(i,j)\end{displaystyle}$, the morphism from $Hom_{\begin{displaystyle}2\varinjlim2\varprojlim\end{displaystyle}}(\mathcal{X},\mathcal{Y})$ to $Hom_{\begin{displaystyle}2\varinjlim2\varprojlim\end{displaystyle}}(\varepsilon(\mathcal{X}),\varepsilon(\mathcal{Y}))$ induced by $\varepsilon$ is the natural morphism :
$$\begin{displaystyle} \varinjlim_{i \in \Ic}\varprojlim_{j \in \Jc}Hom\big(\cat a(s, Id_{j})X_{j
}, \cat a(s', Id_{j})Y_{j}\big)\end{displaystyle} \longrightarrow \begin{displaystyle}\varprojlim_{j \in \Jc} \varinjlim_{i \in \Ic}Hom\big(\cat a(s, Id_{j})X_{j
}, \cat a(s', Id_{j})Y_{j}\big)\end{displaystyle}$$
As filtered colimits commute with finite limits, the morphism above is an isomorphism and $\varepsilon$ is fully faithful. \\

The proof that $\Psi$ is essentially surjective is similar to the proof of the commutation of filtered limits and finite colimits. \\
Let $\Big(\{i_{j},X_{j}\}, \{[\theta^\mathcal{X}_{t}]\}\Big)$ an object of $\begin{displaystyle}2\varprojlim_{j \in \Jc}2 \varinjlim_{i \in \Ic}\aaa(i,j)\end{displaystyle}$. \\
Using  the property $(ii)$ of a filtered category inductively,  one proves that there exist  an object $k' \in \Ic$ and, for any $j \in \Jc$, a morphism $s'_{j} : i_{j} \rightarrow k'$ in $\Ic$.
Thus, for all $t$ morphism of $\Jc$, $[\theta_{t}^\mathcal{X}]$ can be viewed as the class of an object $\vartheta_{t}^\mathcal{X}$ of 
$ Hom\big(\cat a(s_{j}, Id_{j})X_{i_{j}}, \cat a (s_{j'}, t)X_{j'}\big)$.
Remark that even if the class $[\theta_{t}^\mathcal{X}]$ satisfies the equalities $(1)$ and $(2)$, the objects $\vartheta_{t}^\mathcal{X}$ may not satisfy then. Using the property $(iii)$ of a filtered category inductively, one proves that there exist an object $k$ of $\Ic$ and a morphism $s_{k} : k'\rightarrow k$ such that all the equalities hold also for the objects $\vartheta_{t}^\mathcal{X}$. Hence the triple :
$$\Big( k,\big\{\cat a(s_{k}\circ s_{j}, Id_{j})X_{i_{j}}\big\}, \big\{a_{s_{k},s_{j}}^{-1}\circ \vartheta^\mathcal{X}_{t}\circ a_{s_{k}, s_{j'}}\big\}\Big)$$ 
is an object of $\begin{displaystyle}2 \varinjlim_{i \in \Ic}2\varprojlim_{j \in \Jc}\aaa(i,j)\end{displaystyle}$. For $j \in \Jc$, the objects $X_{i_{j}}$ and $\cat a(s_{k}\circ s_{j}, Id_{j})X_{i_{j}}$ are isomorphic and, for all  morphism $t$ of $\Jc$, we have 
$$[\vartheta_{t}^\mathcal{X}]=[\theta_{t}^\mathcal{X}]$$
Hence we have shown that $\Psi$ is an equivalence of categories.

\cqfd

\bibliographystyle{plain}
\bibliography{biblio}

\end{document}